\documentclass[times]{qjrms4}

\usepackage[colorlinks,bookmarksopen,bookmarksnumbered,citecolor=red,urlcolor=red]{hyperref}
\usepackage{mathtools}
\newcommand\BibTeX{{\rmfamily B\kern-.05em \textsc{i\kern-.025em b}\kern-.08em
T\kern-.1667em\lower.7ex\hbox{E}\kern-.125emX}}

\usepackage{moreverb}

\def\volumeyear{2013}

\begin{document}

\runningheads{J.~M.~Restrepo}{Dynamic Likelihood Filter}

\title{The Dynamic Likelihood Filter}

\author{J.~M.~Restrepo\corrauth}

\address{Department of Mathematics \& College of Earth Ocean and Atmospheric Sciences, 
Oregon State University}

\corraddr{Department of Mathematics, Oregon State University, Corvallis OR 97330}

\begin{abstract}
A  Bayesian data assimilation scheme is formulated for advection-dominated or hyperbolic  evolutionary problems, and observations.
The method is referred to as the dynamic likelihood filter because it exploits the model physics to dynamically update 
the likelihood with the aim of making better use of  low uncertainty sparse  observations.

The filter is applied to  a problem with linear dynamics and Gaussian statistics, and compared to the exact estimate, a model outcome, and the Kalman filter estimate. Its estimates are shown to be superior to the model outcomes and the Kalman estimate, when the
observation system is sparse. The added computational expense of the method is linear in the number of observations and thus
computationally efficient, suggesting that the method is practical even if the space dimensions of the physical problem are
large.

\end{abstract}


\maketitle


\section{Introduction}

Data assimilation is a Bayesian framework for combining models and observations, taking into account their respective inherent uncertainties.  These uncertainties are often represented by statistical processes. Among other things,  model uncertainties might represent unresolved, or poorly understood physics, or poorly constrained parameters, or uncertainties in boundary conditions and initial conditions. Observation uncertainties derive from  measurement errors and errors associated with projecting the observations onto 
model space.  The goal of data assimilation is to produce estimates of moments of the posterior distribution of the state vector, 
conditioned on  observations. 

Whether the problem is time-dependent or not, under linear dynamics and when the noise processes  are Gaussian, 
the minimizer of the posterior distribution covariance can be found via least squares. 
For time-independent problems an efficient computational optimization strategy, based upon variational extremization,  is known as 3DVAR  (see \cite{hamillsnyder}).
For time dependent problems a {\it smoother} estimate can be found via 4DVAR  (see \cite{Court93} , and \cite{LorHam} and references contained therein), or via the sequential Kalman smoother (see \cite{wunschbook}).
 
 Geophysical processes are seldom linear. There are data assimilation methods that can handle nonlinear/non-Gaussian
problems, however, none seems to be capable of handling the inherently large  number of degrees of freedom of evolution problems that can be captured by partial differential equations at reasonable resolutions, {\it e.g.}, a weather model. 
Among the methods capable of handling
the nonlinear/non-Gaussian cases  we could mention sampling-based methods, such as the Path Integral Monte Carlo 
of  \cite{path} and  \cite{drifter}, and the Langevin sampler ({\it cf.}, \cite{voss}); a host of particle filter methods ({\it cf.},  \cite{KERAJ},  and \cite{implicitparticle} and references contained therein); variance-minimizing methods (see \cite{kush1}, \cite{kush2}, \cite{kush3}); entropic methods, such as the 
 Mean Field Variational Method of  \cite{ERA-I}.  Ad-hoc approaches that yield reasonable outcomes in  high dimensional weakly nonlinear problems, with strong statistical unimodality,  are also known. Of note is  the ensemble Kalman Filter (enKF, see \cite{eve,ev04}). 
  Its ad-hoc nature rests upon the fact that the method uses nonlinear dynamics in the forecast, and Gaussian assumptions in the analysis, {\it i.e.}, the stage in which data and model outcomes are blended. The method exploits ensemble ideas in order to 
 estimate covariances, hence, it is very efficient in computational storage requirements.
The method is also attractive because it easily handles dynamic problems that present themselves 
  in the form of legacy code.

Filtering techniques have been proposed that make use of  some form of nudging  to improve the stability of the filtering method in the presence of nonlinearity, and/or when there are large discrepancies in the relative uncertainty in model outcomes and measurements. 
The earliest strategies achieve nudging via empirical dynamic relaxation (see \cite{hoke76}, for example). More modern strategies that are worthy of mention: the equivalent-weights particle filter \cite{mades}, and the near-optimal guided particle filter  of \cite{weare13}. A different approach,  a predictor-corrector method, was proposed recently by  \cite{disp}, is called
 {\it displacement assimilation}. In this method  kinematic flow constraints are used in the filtering of data and model outcomes. Specifically, the estimation process is constrained  using the divergence-free condition, a fundamental property of the velocity field underlying the advection/dispersion of a tracer in an incompressible flow. The incompressibility condition is imposed on the filtering forecast of the time-dependent problem via a correction to the local metric of space/time and thus adds to the forecast and analysis stage, another stage that makes a displacement correction. Similar methodologies have been proposed, for example, by Ravela and collaborators ({\it cf.}, \cite{ravela12}),    and \cite{nwha}.

A data assimilation strategy based on  filtering is proposed  that is tailored to evolutionary problems dominated by 
advective processes, in situations where  low uncertainty but sparse observations are available. Sparsity in the observation
system is the norm, rather than the exception, in spatially extended, temporal forecasting and analyses. 
The proposed method is called the {\it dynamic likelihood filter} (DLF) because it uses the
dynamics of the physical problem to update in time the likelihood probability density function. By doing so a sparse network of low
uncertainty observations will be ``nudged" in order to have a greater impact on the posterior distribution of the state vector, conditioned on observations. The DLF  can be used in forecasting, since it can project present and past observations into the future in order to 
 retain the capability of making Bayesian estimates.
Two versions of the methodology are foreseen: the case where the observation system is stationary, and when it has its own dynamic.
The stationary case will be addressed in this study, applied to the simplest of physical models. The moving observation system, to 
be addressed separately, leads to a variant of an assimilation scheme called Lagrangian data assimilation (see  \cite{kij1}, for
example, \cite{drifter} and \cite{DKF},  for issues addressing non-Gaussianity).

\section{Statement of the Problem}

We consider the problem of generating moment histories, over a time spanning $t \in [0, t_f]$,  of the ${\mathbb R}^N$  random state vector $U(t)$, 
conditioned on  noisy observations $Y(t_m) \in \mathbb{R}^K$, $m=1,2,...,M$.  The time $t_f$, for the DLF, may be set in the future.
The measurement or filtering times $t_m$ do not extend beyond the present time $t_{p}$.
The simplest possible linear dynamics, with uncertainties represented by normally-distributed  noise, will be pursued here.

For notational convenience the Kroneker delta function will be redefined as follows: 
\[
\delta_{n,m}:= \left\{ \begin{array}{l} 1, \quad \mbox{if} \, \, t_n=t_m, \\
                                                 0, \quad \mbox{otherwise.}
                                                 \end{array} \right. 
\]

 \subsection{The Dynamics}
 
It is assumed that   $U(t)$ is a finite-dimensional vector, with entries $U_\ell(t) = u(x_\ell,t)$, 
 $\ell = 1, 2, ..., N$. The coordinates $\{x_\ell \}_{\ell=1}^N:=X$ are fixed, distinct, and equi-distant (though they do not have to be so). The $X$ coordinates will be denoted the {\it state coordinates}.
 The wave $u(x,t)$  obeys
\begin{eqnarray}
u_t - C(x,t) u_x &=& F(x,t), \quad t>0, x \in[0,L], \nonumber \\
\langle u(0,t) \rangle &=&  \langle u(L,t) \rangle,  t \ge 0, \nonumber \\
u(x,0) &=& {\cal U}(x), \quad  x \in[0,L],
\label{waveeq}
\end{eqnarray}
 the subscripts $x$ and $t$ connote partial differentiation with respect to these variables, the  
 $<\cdot>$ is the ensemble-average operator.
It will be further assumed that $ \langle F \rangle$ and $ \langle {\cal U} \rangle$  are periodic on the interval $L$. The  source/sink or {\it forcing} $F$ as well as the phase speed $C(x,t)$ are noisy, with  known, zero-mean, additive noise. 

The solution of (\ref{waveeq}) for a finite set of solutions following characteristics can be found as follows:
Let the vector $\Phi(t)$ be such that $\Phi_\ell(0)={\cal U}(x_\ell)$. For $\ell=1,2,...,N$,
\begin{eqnarray}
\frac{d \Phi_\ell}{dt}&=&F_\ell(x,t),  \quad t> 0, \nonumber \\
\Phi_\ell(0) &=& {\cal U}(x_\ell).
\label{eq1}
\end{eqnarray} 
 The initial condition ${\cal U} (x_\ell)$ has a known probability density function (pdf), and
 \[
 F_\ell \, dt =f(x,t) dt + A_\ell(t) dW^{(f)}_\ell(t),
 \]
  where $\langle F_\ell(x,t) \rangle =f(x,t).$
The noise in the forcing is captured by the incremental zero-mean Wiener process $dW^{(f)} (t) \in \mathbb{R}^N$, with known variance $A(t)  \in \mathbb{R}^N$. The solution of (\ref{eq1}) requires the solution of the  equation for the characteristics, namely,
\begin{eqnarray}
\frac{dx_\ell(t)}{dt} &=& C_\ell(x,t), \quad t> 0, \nonumber \\
x_\ell(0) &=& x_\ell, \quad \ell = 1,2,...,N,
\label{characteristics}
\end{eqnarray}
where 
\[
C_\ell(x,t) \, dt = c(x,t) dt + B_\ell(t) dW^{(c)}_\ell,
\]
with
$c(x,t):=\langle C \rangle(x,t)$.
The second term in the characteristic equation represents uncertainties in the wave speed. These uncertainties are 
assumed to be well captured by the (zero mean) Wiener incremental process $dW^{(c)}(t) \in \mathbb{R}^N$, with known variance $B(t) \in \mathbb{R}^N$.

The solution of (\ref{waveeq})-(\ref{characteristics}) will be referred to as the  {\it Exact} solution. 
 Ensemble   members we
denote as {\it Truth} will be constructed, using the exact solutions.
 The Truth estimate  is an ensemble member of the exact solution set, linearly interpolated onto the grid $X$, 
 evaluated at discrete times $t_n=n \, \Delta t$, $n=0,1,2,.., N_f$, where $\Delta t$ is constant.

In geoscience applications, evolutionary models for physical phenomena commonly present themselves in the form of a
computer code, for the discrete-in-time and -in-space approximation of a model for the physics in question. 
The  {\it Model} outcome $V_n \in \mathbb{R}^N$ will be a discretized approximation of the vector $U(t)$ representing
approximations of the solution on $X$ of (\ref{waveeq}) at times
$t_n$. 
The discretization generates 
truncation errors which manifest themselves as biasing errors. This source of error  will not be addressed here. 
There are errors and uncertainties associated 
with how well the unresolved physics of the process is captured. In hyperbolic problems it is typical to have uncertainties in 
the forcing {\it i.e.}, sources/sinks,  as well as in the wave speed.

 The model is assumed to be of the form
\begin{eqnarray}
\sqrt{\Delta t} \Delta w &=& -V_{n} + L_n V_{n-1} + \Delta t f_{n-1},  \quad n=1,2,\ldots, N_f, \nonumber \\
\label{eq3}
\end{eqnarray}
where $V_n$ is periodic, {\it e.g.},  $\langle V_n^{0} \rangle = \langle V_n^{N} \rangle$, and the periodic
 initial condition for the vector is known: $V^\ell_0:={\cal U}(x_\ell)$. Samples of the incremental noise  $\Delta w$ come from a normal variate (the prior  $\pi_n$ is a normal variate). The model noise  variance $Q_n =  \Delta t \langle w_n w_{n'}^{\top} \rangle  \delta _{n,n'}$ is assumed known (that the noise is uncorrelated is not a requirement in what follows).

At discrete times,  
 {\it Observations} are available:
\begin{equation}
 \epsilon(t_m)  = Y(t_m)-H(t_{m}) V(t_{m}), \quad m=1,...,M.
\label{data}
\end{equation} 
The observation errors 
$\epsilon(t_m)$  have a  known pdf. The distribution of $\epsilon_m$ is given by the likelihood $\pi_m(Y_m|V_m)$.
 A linear observation network is considered. The observation network will be
  assumed to coincide with  a subset of or with all locations specified by  $X$.  The  observation matrices, $H(t_m): 
  {\mathbb R}^N  \rightarrow {\mathbb R}^K$. 
 Without loss of generality it will be assumed that observations are taken at discrete  intervals $t_m \in[0,t_p]$. Further, the observation  or {\it filtering times} are taken to be equally spaced, $t_m = m \, \delta t$,   $m=1,2,..,M$, and $\delta t \ge \Delta t$, fixed. 
It is assumed that the model and observation noise processes are uncorrelated.
The observation errors are normally distributed and have variance 
\begin{equation}
 R_m:=\langle \epsilon_m \epsilon_{m'}^\top \rangle \delta_{m,m'}.
 \label{meascov}
 \end{equation}


\section{Data Assimilation via Filtering}
\label{da}
The data assimilation problem can be stated as follows: {\it  Find an estimate of the ensemble mean 
and the uncertainty of the random
vector $U(t)$, at times $t_n$, $n=0,1,...,N_f$,  given a set of random observations $Y(t_m)$, $m=1,2,..,M$.}

The estimate will be found by minimizing the trace of  the posterior covariance of the conditional
probability
\begin{equation}
\pi(V|Y):= \pi(Y|V)\prod_{n=0}^{N_f}   \pi_n. 
\label{postkf}
\end{equation}
The prior,  $\pi_n$, informed by model outcomes, depends on $V_i$, $i=n, n-1,...,0$, and is independent of
the likelihood. The likelihood, informed by observations,  is
 \begin{equation}
 \pi(Y|V):=\prod_{m=1}^{M} \pi_m(Y_m|V_m). 
 \label{likekf}
 \end{equation}
The pdfs  $\pi_n$ and $\pi(Y|V)$ are known and will be described below.

 For linear/Gaussian problems obtaining the posterior mean and covariance allows us to fully characterize the posterior distribution. A well-known algorithm  for finding  the Kalman smoother estimator is the RTS algorithm (the algorithm is described in \cite{wunschbook}). It is based on the Kalman Filter  (see \cite{Jazw70}).

Sparsity of observations,  {\it i.e.}, $K \ll N$,  and low uncertainties in them, as compared to the large dimensions and higher modeling uncertainties of   phenomena captured  by partial differential equations, is a realistic situation, particularly in
 geoscience applications. 
  The Kalman filter will be modified in order to improve the estimate of the ensemble-averaged approximation of the wave solution, using  observations made with low uncertainties by a sparsely-distributed measuring system. 
 Sparsity in the observation network can lead
 to a variety of special challenges:
 Sparsity  can lead to likelihoods that are not very informative: the posterior will be overwhelmed by the prior;  if the likelihoods
 are extremely localized, as they would be if the measurements have low uncertainties, it is possible that the Bayesian inference
 problem becomes ill-posed, particularly if there are biases that are not properly accounted for.
Obtaining improved estimates when observations are sparse is the central motivation for formulating the proposed assimilation scheme.  
\subsection{The Kalman Filter (KF)}
\label{kf}

The Kalman Filter produces a sequential estimate in two steps. In the {\it forecast} step the model is used to produce an initial estimate. 
Since the model is linear and the noise is (unbiased) normal, the mean state $\langle V \rangle_{n}$ is obtained from $\langle V\rangle_{n-1}$ using
\begin{eqnarray}
\tilde V  = L_{n-1}  \langle V \rangle_{n-1}  + \Delta t f_{n-1},  \quad n=1,2,\ldots, N_f.
\label{forecast1}
\end{eqnarray}
 It is easy to form the equation for the evolution of the covariance $P_n = \langle e_n e_n^T \rangle$, where $e_n:=V_n - \langle V \rangle_n$: 
 \begin{eqnarray}
\tilde P  = L_{n-1} P_{n-1} L^\top_{n-1}   + Q_{n-1},  \quad n=1,2,\ldots, N_f.
\label{forecast2}
\end{eqnarray}

If no observations are available at time $n$, the posterior is not affected by the likelihood at $t_n$ and thus $\langle V \rangle_n = \tilde V $, 
 and $P_n = \tilde P$. ( $\langle V \rangle_0$, and $P_0$ are known). If, on the other hand, $t_n=t_m$, observations are available. An {\it analysis} step is performed, that takes in the tilde variables $\tilde V$ and $\tilde P$ and updates these to produce a  mean and the covariance estimate. Defined for any step $t_n$,  the analysis step consists of the update 
 \begin{eqnarray}
\langle V \rangle_{n}  &=&  \tilde V   + K_m \left(Y_m - H_m \tilde V \right),  \label{a1} \\
P_n &=& (I_N- K_m H_m )\tilde P.
\label{analysis}
\end{eqnarray}
 The second term in (\ref{a1}) is called the {\it innovation}.
In (\ref{analysis})  $I_N$ is the $N$-dimensional identity matrix. The Kalman Gain is defined as
 \begin{equation}
 K_m = \tilde P H^\top_m \left[ H_m \tilde P H_m^{\top} + R_m \right]^{-1} \delta_{n,m}.
 \label{kalman}
 \end{equation}

\subsection{The Dynamic Likelihood Filter (DLF)}
\label{dlkf}

An alternative Bayesian statement is proposed, for the posterior pdf. In the DLF the pdf, at time $n$,  is
\begin{equation}
 \pi(Y|V)_n:=\pi_m(\{{\cal H}Y_m\}_{t_m \le t_n}|V_n),
 \label{likedlf}
\end{equation}
where the data $Y_m$, for all measuring times $t_m$ less than or equal to time  $t_n$, can influence the likelihood at time $t_n$
({\it cf., } to (\ref{likekf})).  The linear operator ${\cal H}$ projects data from the past to the present and further, maps it
onto the state coordinates, $X$.

In other words,  the observations will not only inform the likelihood at time $t_m$, but will do so at subsequent times.
The observation data and its uncertainty are propagated forward in time, from $t_m$ to $t_n$. The forward-propagated data 
is used, provided that their  inherent uncertainty has not grown (degraded) beyond some set threshold. In doing so
high quality but sparse resources of information are exploited more thoroughly.  There is good reason to think that in  highly local problems in which information travels at finite speeds, this strategy may deliver improved estimates, particularly of phase-sensitive information as well as structure. Clearly, if the data is neither sparse or is endowed with  high uncertainties, compared to the
uncertainties in the model, the methodology would not be recommended.

One can envision two ways of using the likelihood dynamically: one is when the observation network is stationary. 
The other modality would correspond to when the observation network moves in space and 
time. This latter case corresponds to  a type of  Lagrangian data assimilation ({\it cf.}, \cite{kij1}), and will be developed in detail
in a separate paper.

\subsubsection{Dynamic Likelihood on  a Fixed Observation System}

What is meant by a stationary observation network is that the observations take place at prescribed locations
$Z \subseteq X$, for all $t$, where dim$(Z) = K$. 
The {\it observation coordinates} $Z$  are taken to coincide with components of the lattice $X$ in order to avoid the complications
that come with the addition of a Gauss-Markov interpolation between the state coordinates and the observation coordinates.
(A moving observation system, on the other hand, will have $Z=Z(t)$). 

Observations from the fixed observation system will be used at times $t_m$ when they become available, 
as well as at subsequent times, once their phase and uncertainty is updated. Using the known estimate for the mean
speed $c(x,t)= \langle C(x,t) \rangle$, the observation phase is updated by computing the wave characteristics
emanating from the location of the observation.
The semi-Lagrangian approximation for the virtual dynamics of observational data is used 
 to propagate the observations into future fixed times, the method makes use of the  semi-Lagrangian approximation  \cite{staniforth} to solutions to data that evolves according to $Y_t - c(x,t) Y_x = 0$, on the time grid. That is,
\begin{eqnarray}
\zeta_{n+1} &=& \Delta t c(\zeta_n,t_n) + \zeta_n,  \quad t_n \ge t_m, \nonumber \\
Y(\zeta_{n+1},t_{n+1}) &=& Y (\zeta_{n},t_{n}),
\label{semilagrangian}
\end{eqnarray}
with $\zeta_0 =H(t_m) X$,  and $Y(\zeta_0,t_m)=Y_m$.  

The variance associated with the measurements at $t_m$  is known and equal to $R_m$, but it cannot be expected to remain the same, for $t_n>t_m$.
Since the inherent dynamics of the problem are being used to project forward in time, it may be argued that the manner in which the uncertainty of these projected observations change is dictated by the uncertainty in the dynamics. If this argument is reasonable
\begin{equation}
R^{n+1}_m =   A_n(t) [A_n(t)]^\top \Delta t + R^n, \quad t_n \ge t_m,
\label{Req}
\end{equation}
where $R^m_m = R_m$.

The decision on how this data is used must rely on sensible choices in the computation and the estimation process.
It is decided that measurements propagated forward in time that reach some threshold of uncertainty should cease contributing to the likelihood.  Presuming the data has lower uncertainty than the model, a reasonable uncertainty threshold for measurements propagated forward in time, is that their variance not exceed the variance in the 
model in the same spatial neighborhood. Using (\ref{Req}) it is easy in fact to estimate the total time 
beyond $t_m$ that data is viable, given this threshold or some other criteria.

The dimension of the viable observations can exceed $N$. Computationally, this may lead to having to
address numerical stability issues differently when the system goes from under determined to an overdetermined system. Moreover, the Kalman strategy presumes that a single innovation vector be used to update the forecast
which in turn means addressing complexities in the observation matrix.  In what follows choices are made to 
keep the number of observations used at any given time $t_n$ bounded by $N$. Furthermore, the propagated
observations are ordered by the norm of their uncertainty and preference will be given to observations with lower uncertainty.
 
\subsubsection{When Measurements Relate to the State via a Linear Transformation}
\label{LT}

For the linear dynamics  the Kalman Filter is still a viable estimation strategy.
In the DLF, as applied to linear/Gaussian problems,  the forecast stage of the Kalman Filter remains the same. 
What changes is the analysis stage. The correction stage of the algorithm will be referred to as 
the multi-analysis. It entails building an alternative innovation vector and a Kalman gain.

In the usual Kalman Filter, an ``observation" matrix projects the  state vector onto the space of measurements and their uncertainties. Here, the ``reverse" case is developed, wherein a matrix projects the measurements and their errors onto the state vector.

The crux is to develop the sequential estimation when the observation matrix ${\cal H}$ appears thusly:
\begin{equation}
{\cal H}^n_m Y^n_m = V_n + {\cal H}^n_m  \epsilon^n_m,
\label{thedata}
\end{equation}
at time $t_n \ge t_m$. Here
$\epsilon_m^m$ is equal to $\epsilon_m$.
The multi-analysis stage is now
\begin{equation}
\langle V \rangle_n = \tilde V + {\cal K}_m ({\cal H}_m Y_m - \tilde V) \delta_{m,n}.
\label{dlf2}
\end{equation}
Since $P_n = \langle (\langle V_n - \langle V \rangle_n)(\langle V_n - \langle V \rangle_n)^\top \rangle$, it is then possible to express the covariance,
using (\ref{thedata}), as
\begin{equation}
P_m = \langle (I- {\cal K}_m ) \tilde P  (I- {\cal K}_m )^\top \rangle + {\cal K}_m {\cal H}_m R_m {\cal H}_m^\top {\cal K}_m^\top .
\label{pneq}
\end{equation}
The trace of the covariance is
\[
\mbox{Tr}[P_m] = \mbox{Tr}[\tilde P] - 2 \mbox{Tr} [{\cal K}_m \tilde P] + \mbox{Tr}[{\cal K}_m (\tilde P + {\cal H}_m R_m {\cal H}_m^\top) {\cal K}_m^\top].
\]
Differentiating with respect to ${\cal K}_m$ and setting the derivate to zero, one finds the extremizer of the trace is
\begin{equation}
{\cal K}_m = \tilde P (\tilde P + {\cal  H}_m R_m {\cal H}_m^\top)^{-1} \delta_{m,n}.
\label{TheK}
\end{equation}
Using (\ref{TheK}) back in (\ref{pneq}) we find the update to the covariance is
\begin{equation}
P_n = (I -  \delta_{m,n} {\cal K}_m) \tilde P.
\label{updatep}
\end{equation}

\subsection{Formulating a Low Uncertainty Likelihood}
The expressions for the Kalman Gain (\ref{TheK}) and the uncertainty (\ref{updatep}), for this linear-Gaussian version of the DLF,
 are purely formal at this point. 
 In the dynamic likelihood assimilation it will be generally the case that several observations will have a bearing on a single state variable entry. If the measurements are very sparse in space, a simple option is to use 
\[
{\cal H}_m^n Y_m^n (x_\ell) := {\cal I} {\cal J}
Y^n_{m}(x_j(t_n)) = {\cal I} {\cal Y}(x_\ell),
\]
The {\it projection operator} ${\cal J}$ is defined as 
\begin{equation}
{\cal J}:= [\delta^n_{\ell,j} (1-b^n_j)+ \delta^n_{\ell-1,j} b^n_j], 
 \label{project}
 \end{equation}
where $b^n_j = b_0  \,\mbox{rem}(x_j(t_n),\Delta x)$, the remainder, and $b_0$ is a normalization constant. 
In (\ref{project}) the Kroneker delta symbols are used in the traditional sense. (The projection operator ${\cal J}$
 must be modified to take into account boundary conditions).
The matrices project the updated datum at position $x_j$ onto the fixed grid $x_\ell$ by weighting the data according to the proximity of $x_j$ to the grid on which $x_\ell$ belongs. A simpler alternative is to set $b_j^n$ to zero. This latter option will be used in the example calculations that appear later on.

With regard to choosing which measurements to use in the filtering stage, when there are choices to be made, the following 
procedure is suggested: at any given time, measured or updated measurements are ordered in increasing size in their associated
uncertainty. The rank-ordering operator accomplishes this. 
Observations available presently or propagated from the past may be available to perform a multi-analysis and may share the same 
projected location $x_\ell$.
The  {\it rank-ordering operator}   ${\cal I}$ assigns a single datum to $x_\ell$, among those that are  available. The datum will have
 the lowest uncertainty among all of  the data that can inform the likelihood at location $x_\ell$, from the present or from the past.
 Let ${\cal I}_k = D_{k+1} + E_{k+1}$, $k=0,1,...$, then
 \begin{equation}
 {\cal I}_0 {\cal Y}= D_1 {\cal Y}+ E_1 [ D_2 {\cal Y}+ E_2 [D_3 {\cal Y}+E_3[ D_4 {\cal Y}+..] ...]],
 \label{ranker}
 \end{equation}
{\it i.e.}, the rank-ordering operator is applied recursively. 

\begin{figure}[ht]
 \centering
\scalebox{0.4}{\includegraphics{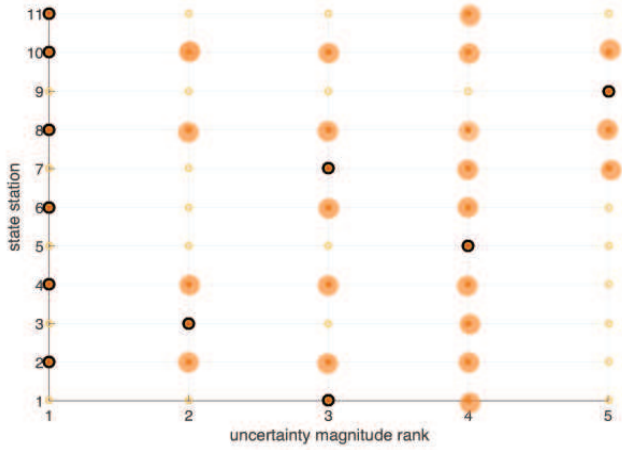}}
\caption{Rank-ordering operator schematic. Measurements, projected forward or from the present, are organized from left to 
right in increasing uncertainty. The measurements have been already been projected onto the state stations, via the ${\cal J}$
operator. The state stations are labeled 1 through 11.  Available data appear as filled circles. The data that will be used in the innovation appears as dark, multi-colored filled circles. The datum with the least uncertainty at a specific state station will contribute to the innovation.
 In this example data manages to inform all of the stations.}
\label{rank}
\end{figure}

An example illustrates the rank-ordering operator.   Suppose there are 11 state stations. The goal is to construct the innovation 
vector at some time step $n$. Figure \ref{rank} schematically portrays the observations and/or their projected values. It will be
assumed that the uncertainty associated with measurements grows linearly in time. The measurements, or their 
dynamic updates, have been 
projected already onto the state stations.  Figure \ref{rank} displays the set of observations and/or their projected values, at time $n$.
The data has been   arranged, from left to right, in order of increasing uncertainty. 
For this example the matrices are 
\begin{eqnarray*}
&& D_1 = \delta_{i,j}
\left[  \begin{array}{c}  0\\ 1\\ 0 \\1\\ 0\\ 1\\ 0\\ 1\\ 0\\ 1\\ 1 \end{array} \right], 
D_2=\delta_{i,j} \left[ \begin{array}{c} 
0\\ 1\\ 1 \\1 \\0 \\0\\ 0\\ 1\\ 0\\ 1\\ 0  \end{array} \right], 
D_3=\delta_{i,j}  \left[ \begin{array}{c} 
1\\ 1\\ 0 \\1 \\0 \\1 \\1 \\1 \\0 \\1 \\0  \end{array} \right], \\
&& D_4 = \delta_{i,j}   \left[ \begin{array}{c} 
1 \\1\\ 1\\ 1 \\1 \\1\\ 1 \\1 \\0 \\1 \\1 \end{array} \right], 
D_5=\delta_{i,j}  \left[ \begin{array}{c} 
0\\ 0\\ 0\\ 0 \\0 \\0 \\1 \\1 \\1\\ 1\\ 0\end{array} \right], \quad i,j=1,2,...,nx.
\end{eqnarray*}

\subsection{The Mutli-Analysis and the Kalman Gain}

The multi-analysis is 
\begin{equation}
\langle V \rangle_n = \tilde V_n +  {\cal K}_{n} \sum_{m' \in m}  ({\cal H}^{n}_{m'} Y^{n}_{m'} - \tilde V_{n}),
\label{dlfi}
\end{equation}
 The data $Y^n_m$ is advanced via (\ref{semilagrangian}), or its ensemble-forced counterpart if the physics has focing.  The uncertainty $R^n_m$ is updated via (\ref{Req}) (or its forced counterpart).

To obtain an expression for the Kalman Gain, we extremize the trace of the posterior covariance, just as was done in Section \ref{LT}. The
result for the Kalman gain and the uncertainty, respectively,  are
\begin{equation}
{\cal K}_{n} = \tilde P_n  (\tilde P_n +  \sum_{m'\in m} {\cal  H}_{m'}^n R^n_{m'} [{\cal H}^n_{m'}]^\top \vartheta_{m',n}  )^{-1} ,
\label{multigain}
\end{equation}
where  $\vartheta_{m',n}$ is the Heaviside function,
and
\begin{equation}
P_n = (I -{\cal K}_{n}  ) \tilde P_n.
\label{multiupdatep}
\end{equation}

  \subsection{The Linear/Gaussian DLF Algorithm}

The algorithm for DLF, for the linear Gaussian case, which allows the use of the 
Kalman filter framework, proceeds as follows: at some state time $t_n >0$, and given that at $t=0$ 
the mean of the state of the system and its uncertainty are known:
\begin{itemize}
\item   The forecast stage is the same as the one employed in  the Kalman Fillter, {\it i.e.},  (\ref{forecast1})-(\ref{forecast2}).
\item The multi-analysis stage: 
\begin{itemize}
\item Observations from times $t_m < t_n$ are propagated forward up to time at $t_n$. The propagation is effected by 
using the mean eigenspeeds, along characteristics. Their uncertainties are propagated forward as well. At measuring times
the uncertainties in the observations are given by the known measurement error.
\item The innovation vector is built  using the  projection and the uncertainty rank operators.
\item The Kalman gain is computed via (\ref{multigain}).
\item The multi-analysis is performed using (\ref{dlfi}), yielding an estimate of the mean state.
\item The uncertainty of the state is updated via (\ref{multiupdatep}).
\end{itemize}
\item The forecast and multi-analysis is repeated till time $t_f$, along the way the measurements are used in formulating the
multi-analysis stage. If the presumption is that the model error is larger than the observation error, one can use an uncertainty
threshold to shed measurements  in the multi-analysis that, due to being propagated forward in time, have an unacceptably high uncertainty as compared to the model uncertainty at that particular time $t_n$.
  \end{itemize}  
  We emphasize again, that $t_f$ can be set into the future, the filtering takes place using the model and time-projected data.
  
\section{Comparing  the Model, the KF, and the  DLF Outcomes}
\label{numerical}

Comparisons will be made between the outcomes from the Lax-Friedrich calculations, which we will call the ``Model" outcomes, and assimilation results from the Kalman Filter (KF) and the Dynamic Likelihood (Kalman) Filter (DLF).  
These are compared, in turn, to  the Truth.

The observations are available at times $t_m$, on all or a portion of the state space stations $X$.  At time $t_m$ the measurements have  fixed  measurement variance $R$. They are the same for the Kalman and the DLF. However, in DLF the observations are used in the innovation vector at time $t_m$ and times thereafter. For $t>t_m$ the observations are updated in space and time as are their uncertainties. Their impact on the estimate at times $t_n$ and stations $X$ is found via the multi-analysis.

The aims of the calculations are  to show that the mean delivered by data assimilation is superior to the 
model outcome;  and to show that when the observation system is sparse, the DLF yields outcomes
that are better than  KF  with regard to phase, uncertainty, and in many instances, in qualitative terms.

The two parameters in the examples will be the frequency of 
spatial sampling $\xi$, and the frequency of temporal sampling $\tau$. The spatial frequency describes the inverse distance between the fixed observation stations. The temporal frequency gives the rate at which new observations are read. The Model does not make 
use of measurements, the KF makes use of observations as they become available, the DLF uses observations as they become available as well as observations that are updated in time and space.

\subsubsection{Example Problems}
The model problems are chosen because they have analytical solutions. The Exact analytical solutions are linearly interpolated on the $X$ to produce a solution we denote as  {\it Truth}.

 {\it Problem I}:
\[
dx = -\alpha x dt + \beta dW,
\]
with known initial conditions (initial probability density).
In this case the relative drift is the constant $-\alpha$ and the variance of the incremental Wiener process $dW$ is
$\beta$. (In physical applications the constant  $\alpha >0$, where $1/\alpha$ is a commonly referred to as the relaxation time). The
 solution  is
\[
x_{n+1} = x_n e^{-\alpha \Delta t} + \left[\frac{\beta^2}{2 \alpha} (1-e^{-2 \alpha \Delta t}) \right]^{1/2} {\cal N}(0,1),
\]
for $n=0,1,...N_f-1$.
 An approximate solution is given by
\[
x_{n+1} = x_n - \alpha x_n \Delta t + \sqrt {\beta^2 \Delta t} {\cal N} (0,1).
\]
The mean and the variance of the solution are, respectively, 
\[
\langle x_{n+1} \rangle = \langle x_n \rangle e^{-\alpha \Delta t}, \quad \mbox{cov}(x_{n+1}) = 
\mbox{cov}(x_n) + \frac{\beta^2}{2 \alpha} (1-e^{-2 \alpha \Delta t}).
\]

{\it Problem II}:
\[
d x = (\alpha_0 + \alpha_1 t^{1/2})dt + \beta dW,
\]
with known initial conditions. Here $\alpha_0$ and $\alpha_1$ are constants. Problem B has 
a  solution
\[
x_{n+1} = x_n + \alpha_0 \Delta t + \frac{2}{3} \alpha_1 \Delta t^{3/2} +  \sqrt{ \beta^2\Delta t}{\cal N} (0,1),
\]
for $n=0,1,...N_f-1$. ${\cal N}(0,1)$ is a normal variate with variance 1. The mean and the covariance of the
solution are, respectively, 
\[
\langle x_{n+1} \rangle = \langle x_n \rangle + \alpha_0 \Delta t + \frac{2}{3} \alpha_1 \Delta t^{3/2}, \quad \mbox{cov}(x_{n+1}) = \mbox{cov}(x_n) + \beta^2 \Delta t.
\]

\subsubsection{The Model}

The Lax-Friedrichs discretization (\cite{iserles}) of (\ref{waveeq}) is adopted as the Model. Clearly, there are better discretizations than the Lax-Friedrichs, but the point is that in practice we have a 
model, which by definition, is the working approximation of the Truth.
 The forcing will be set to zero in the illustrative calculations. 
The model, for  $u(t_n,X):=U^n \approx V^n$ is 
\[
V^{n+1} =  L^n V^n +\sqrt{\Delta t} \Delta w_n, \quad n=0,1,2, ...
\]
with known $V^0 = \langle U^0 \rangle$ and variance, periodic, and
\[
L^n = \frac{1}{2} (1-\lambda^n_\ell) \delta_{\ell,\ell+1} +  \frac{1}{2} (1+\lambda^n_\ell) \delta_{\ell,\ell-1},
\]
for $\ell=2,...,M-1$, with periodized entries in the first and last rows:
\begin{eqnarray*}
L^n_{1,2} &=& \frac{1}{2} (1-\lambda^n_1) \quad L_{1,M} = \frac{1}{2} (1+\lambda^n_1) \\
L^n_{M,1} &=& \frac{1}{2} (1-\lambda^n_1)  \quad L_{M,M-1} = \frac{1}{2} (1+\lambda^n_1)
\end{eqnarray*}
where 
\[
\lambda_\ell^n = \frac{\Delta t}{\Delta x} c(x_\ell, t_n).
\]

The Courant-Friedrichs-Lewey condition (CFL) is set to $0.99$. In the calculations that follow, once $\Delta x$ is determined the $\Delta t$ is determined via the wave speed and the CFL condition.

\subsubsection{Fixed Parameters in the Examples}
The discretization parameters will remain the same for every example shown. The domain has a length of $L=2$, and it is discretized using $N=50$ grid points in $x$, {\it i.e.}, dim$X$=N=50 is the dimension of the state vector. In time there will be $N_f=200$ time steps, for Problem I and $N_f=100$ time steps, for Problem II. 
The initial condition is a pulse
\[
U^0 = \frac{1}{S}[1- 4 (X-x_0)^2] \mathbb{I}_X + E,
\]
where $E$ is a normal variate vector with uncertainty $P_0=0.02$, and $\mathbb{I}_X$ is a characteristic function.  The normalization $S=\int_0^L [1- 4 (X-x_0)^2] \mathbb{I}_X dx$.
For Example I, $x_0 = 1.25$ and for Example II, $x_0 = 1$.
 The forcing is purely stochastic: $f=0$, $A=0.01$. 
 For Problem I: $\alpha_0=0.1$,  $\alpha_1=0.01$.  For Problem II: $\alpha=0.01$. The variance of the observations is $R=0.02$, and the variance of the model $Q=4 R$. The variance on the wave speed $\beta=0.02$.

\subsubsection{Comparisons}

\begin{figure}[ht]
 \centering
 \scalebox{1}{\includegraphics[height=2.5in,width=3in]{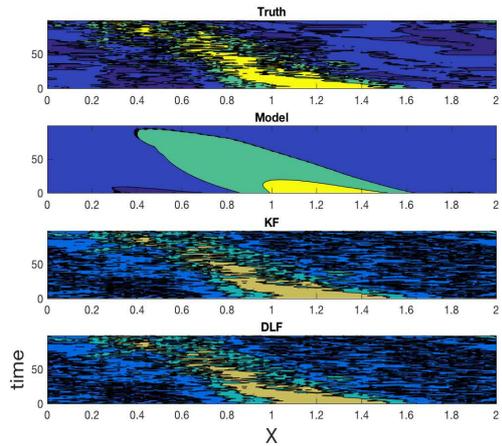}}
\caption{Problem I, Truth, Model, KF and DLF outcomes. New data is read at every time step (the temporal sampling frequency is $\tau=1$) and at every spatial location (the spatial sampling frequency  is $\xi=1$).}
\label{DLKF1a}
\end{figure}

\begin{figure}[ht]
 \centering
 \scalebox{1}{\includegraphics[height=2.5in,width=3in]{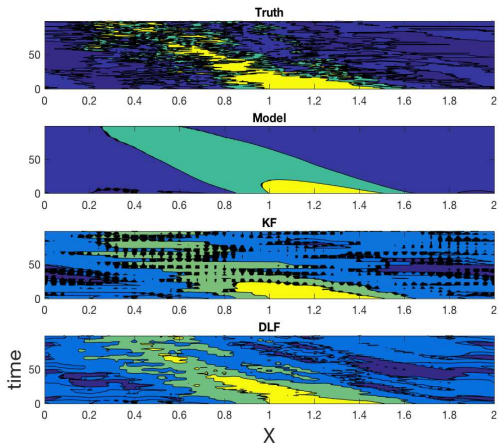}}
\caption{
Problem I, Truth, Model, KF and DLF outcomes; $\tau=1/10$ and $\xi=1/5$.}
\label{DLKF1b}
\end{figure}

\begin{figure}[ht]
 \centering
(a) \scalebox{1}{\includegraphics[height=2.in,width=3in]{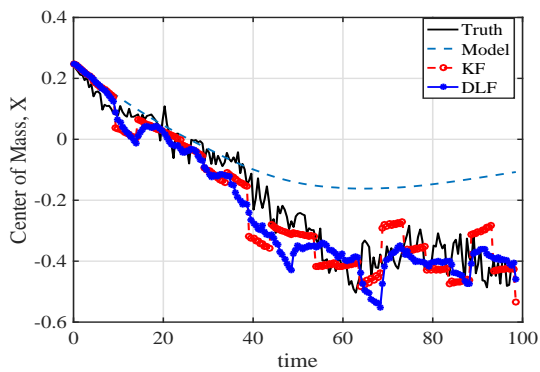}} \\
(b) \scalebox{1}{\includegraphics[height=2.in,width=3in]{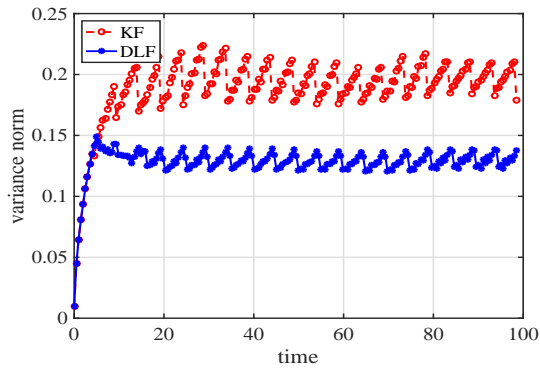}}
\caption{Problem I,  Truth, Model, KF and DLF outcomes, (a) center of mass, (b) uncertainty; $\tau=1/10$ and $\xi=1/5$.}
\label{DLKF1c}
\end{figure}

Figures \ref{DLKF1a}-\ref{DLKF1c} summarize the outcomes for Problem I, with different filter sampling
configurations.  Figure \ref{DLKF1a} shows that the KF and the DLF are indistinguishable, and in very close
agreement with an ensemble member of  the Truth. 

Differences are noticed in the estimates when the filtering is done more sparsely. 
When the sampling is sparse, the DLF is qualitatively similar to the Truth, as compared to the 
KF outcome. This is shown in Figure \ref{DLKF1b}  

In Figure \ref{DLKF1c}a  the evolution of the center of mass of the different outcomes is 
compared, for sparse conditions. 
 The Model estimate, as expected, diverges from the Truth estimate  significantly, however,  the KF or DLF are fairly close to each other and to the mean Truth. It is noted  that, in this example problem over this short period of time, 
  the Model provides a reasonably good estimate of the  phase, when compared to the phase of the Truth. Hence, the KF and
the DLF are given fairly good estimates in the prediction (forecasting) step.  
 The uncertainty in the DLF is smaller than the KF but both of these are well within the expected
 bounds for the posterior uncertainty. See Figure \ref{DLKF1c}b.

\begin{figure}[ht]
 \centering
 \scalebox{1}{\includegraphics[height=2.5in,width=3in]{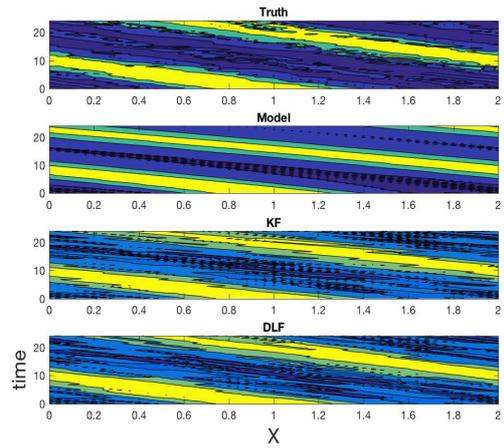}}\\
\caption{Problem II, Truth, Model, KF and DLF mean outcomes; $\tau=1/10$, $\xi=1$.}
\label{DLKF1}
\end{figure}

\begin{figure}[ht]
 \centering
 \scalebox{1}{\includegraphics[height=2.5in,width=3in]{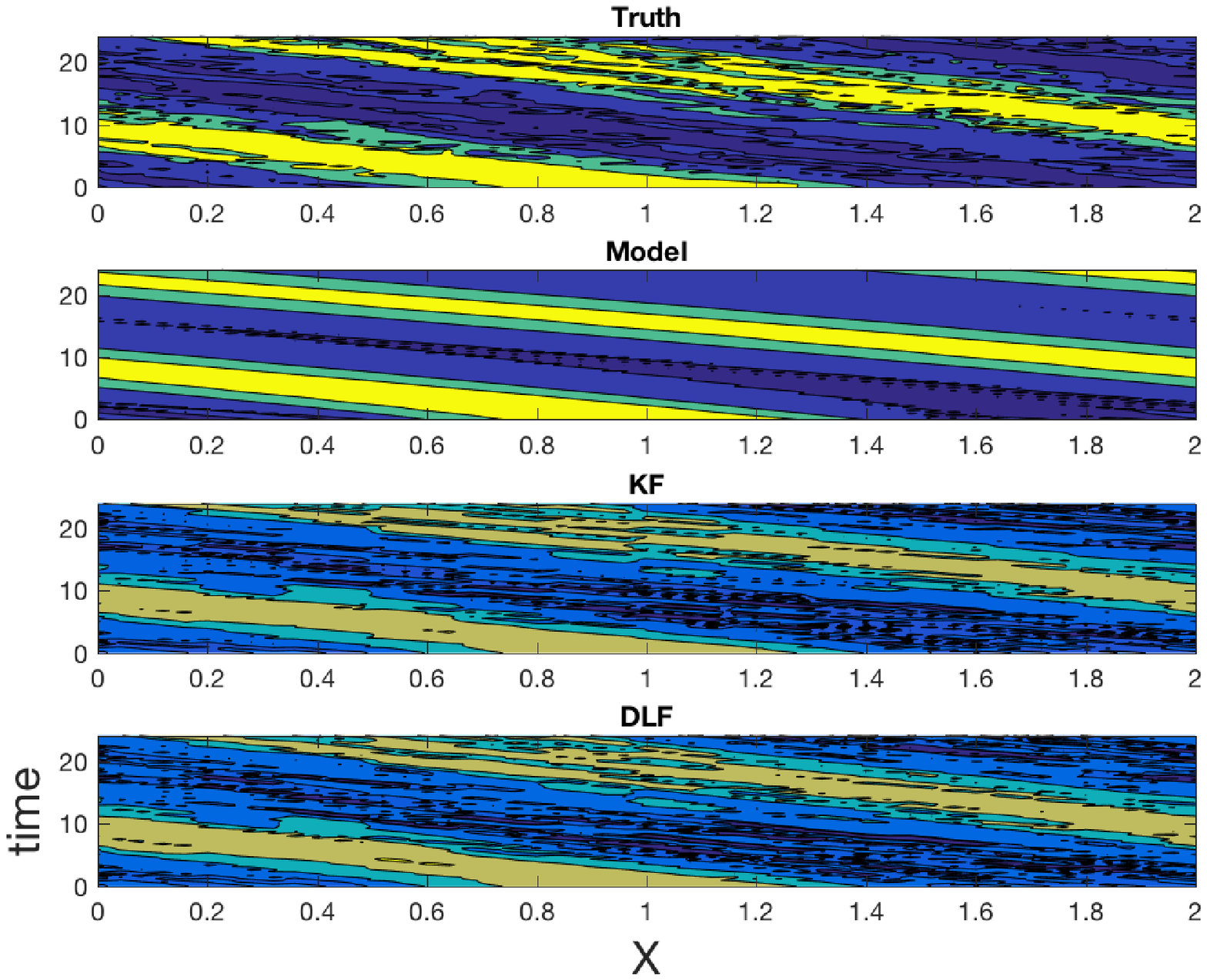}}\\
\caption{
Problem II, Truth, Model, KF and DLF  mean outcomes; $\tau=1$, $\xi=1/4$.}	
\label{DLKF2}
\end{figure}

\begin{figure}[ht]
 \centering
 \scalebox{1}{\includegraphics[height=2.5in,width=3in]{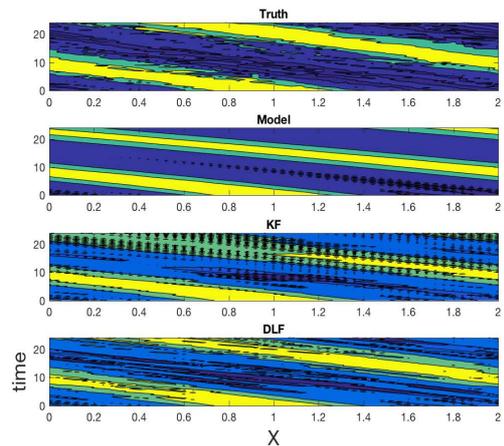}}
\caption{
Problem II, Truth, Model, KF and DLF mean  outcomes; $\tau=1/10$, $\xi=1/4$.}
\label{DLKF3}
\end{figure}

\begin{figure}[ht]
 \centering
\scalebox{1}{\includegraphics[height=2.5in,width=3in]{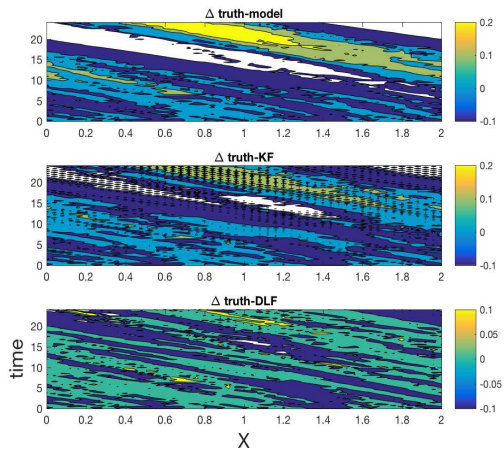}}\\
\caption{Problem II, difference between the Truth and the Model, the Truth and KF, and Truth and DLF; 
$\tau=1/10$, $\xi=1/4$.}
\label{DLKF4}
\end{figure}

\begin{figure}[ht]
 \centering
(a) \scalebox{1}{\includegraphics[height=1.2in,width=2.7in]{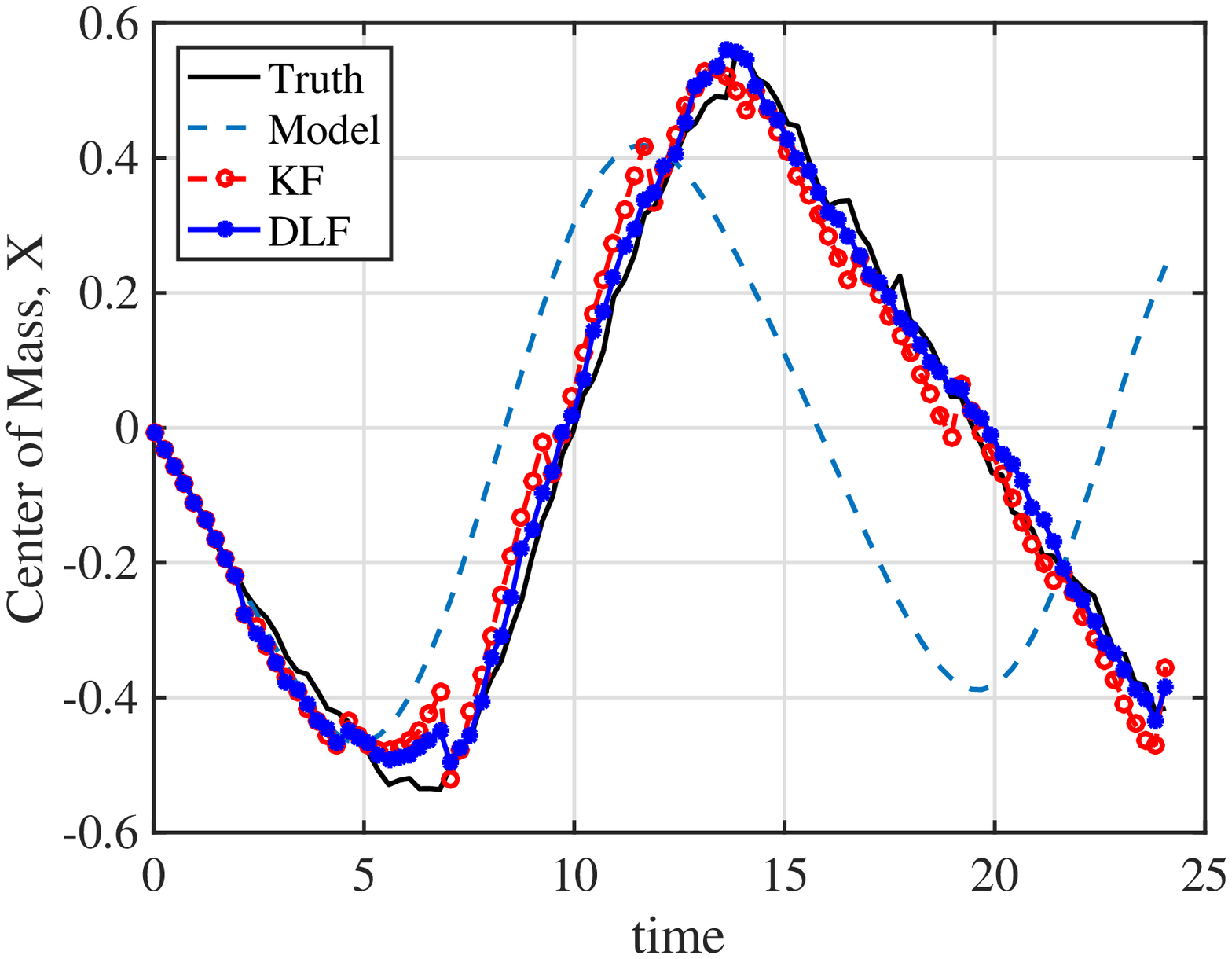}}\\
(b) \scalebox{1}{\includegraphics[height=1.2in,width=2.7in]{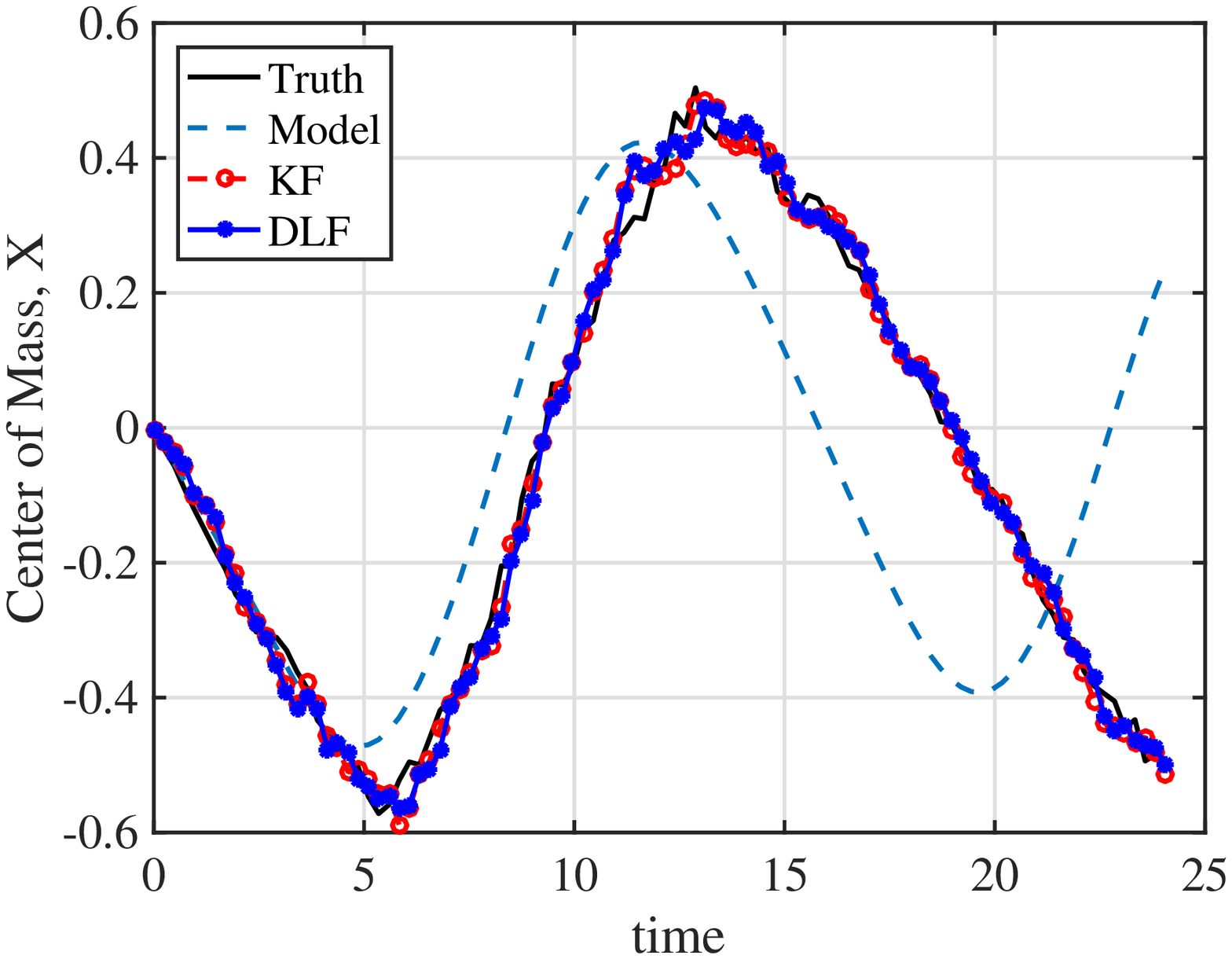}}\\
(c) \scalebox{1}{\includegraphics[height=1.2in,width=2.7in]{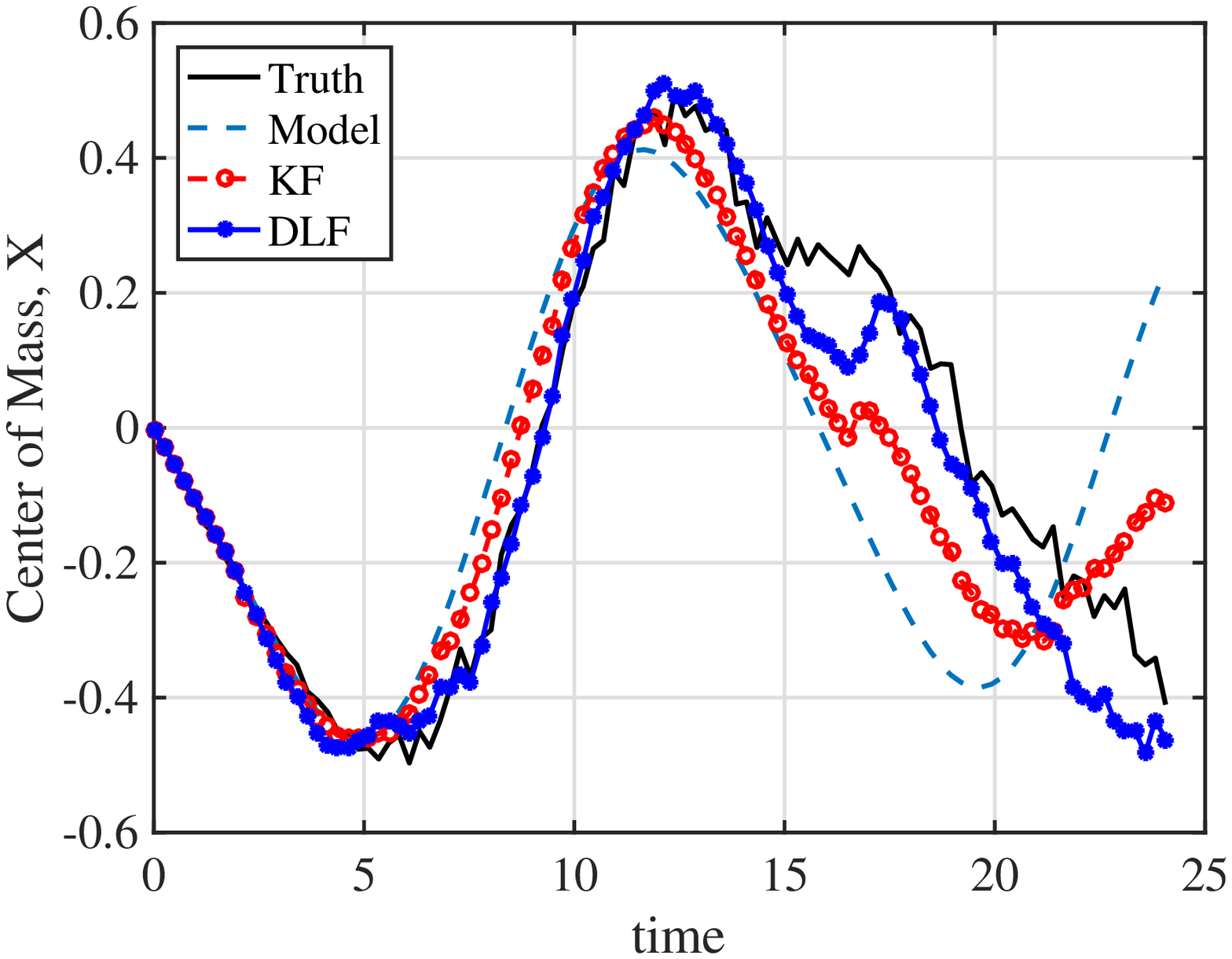}}
\caption{Problem II. Time evolution of the center of mass of the Truth,   
the Model, and the KF, and DLF  estimates;   (a) $\tau=1/10$, $\xi=1$, (b) 
$\tau=1$, $\xi=1/4$, and (c) $\tau=1/10$, $\xi=1/4$.}
\label{DLKFphase}
\end{figure}

\begin{figure}[ht]
 \centering
(a) \scalebox{1}{\includegraphics[height=1.2in,width=2.7in]{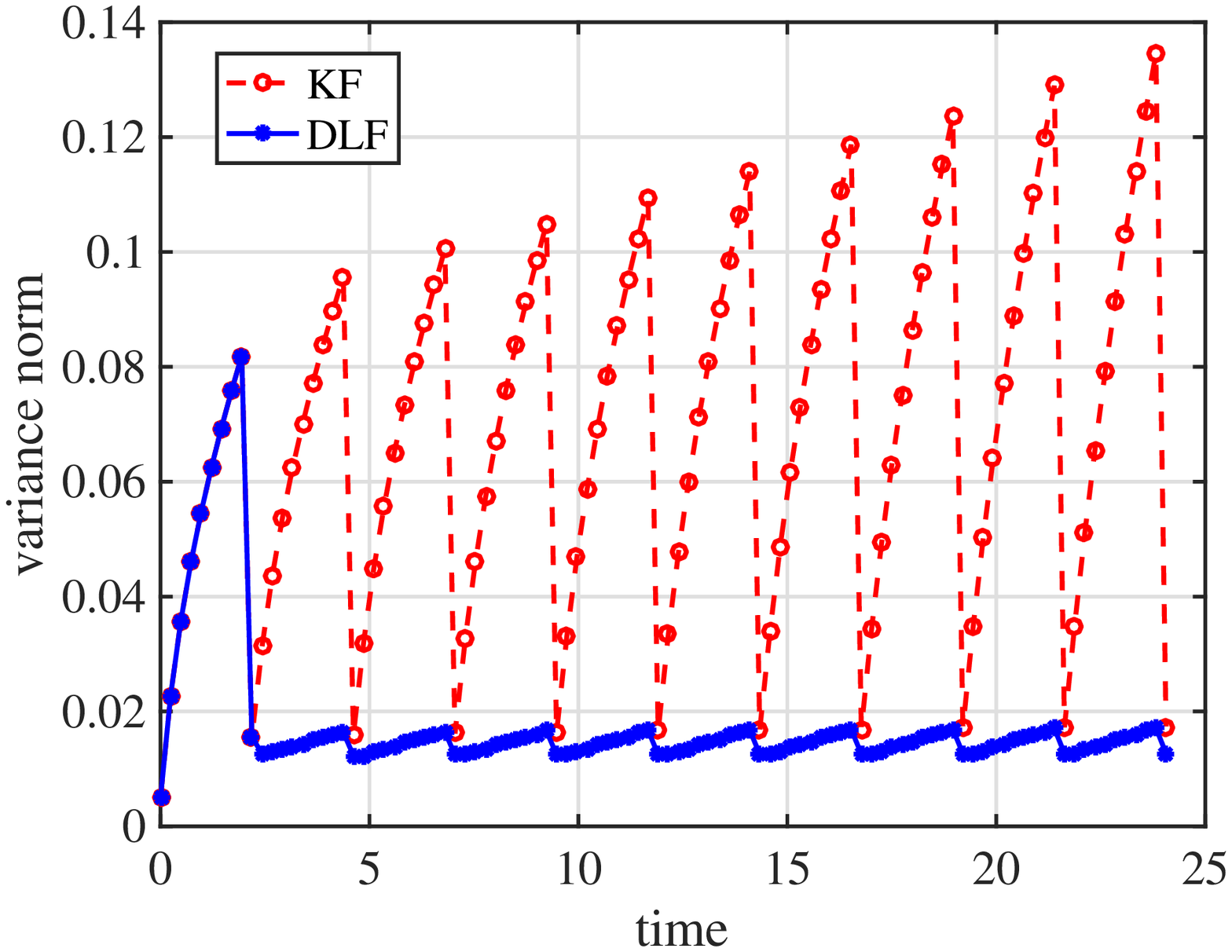}}\\
(b) \scalebox{1}{\includegraphics[height=1.2in,width=2.7in]{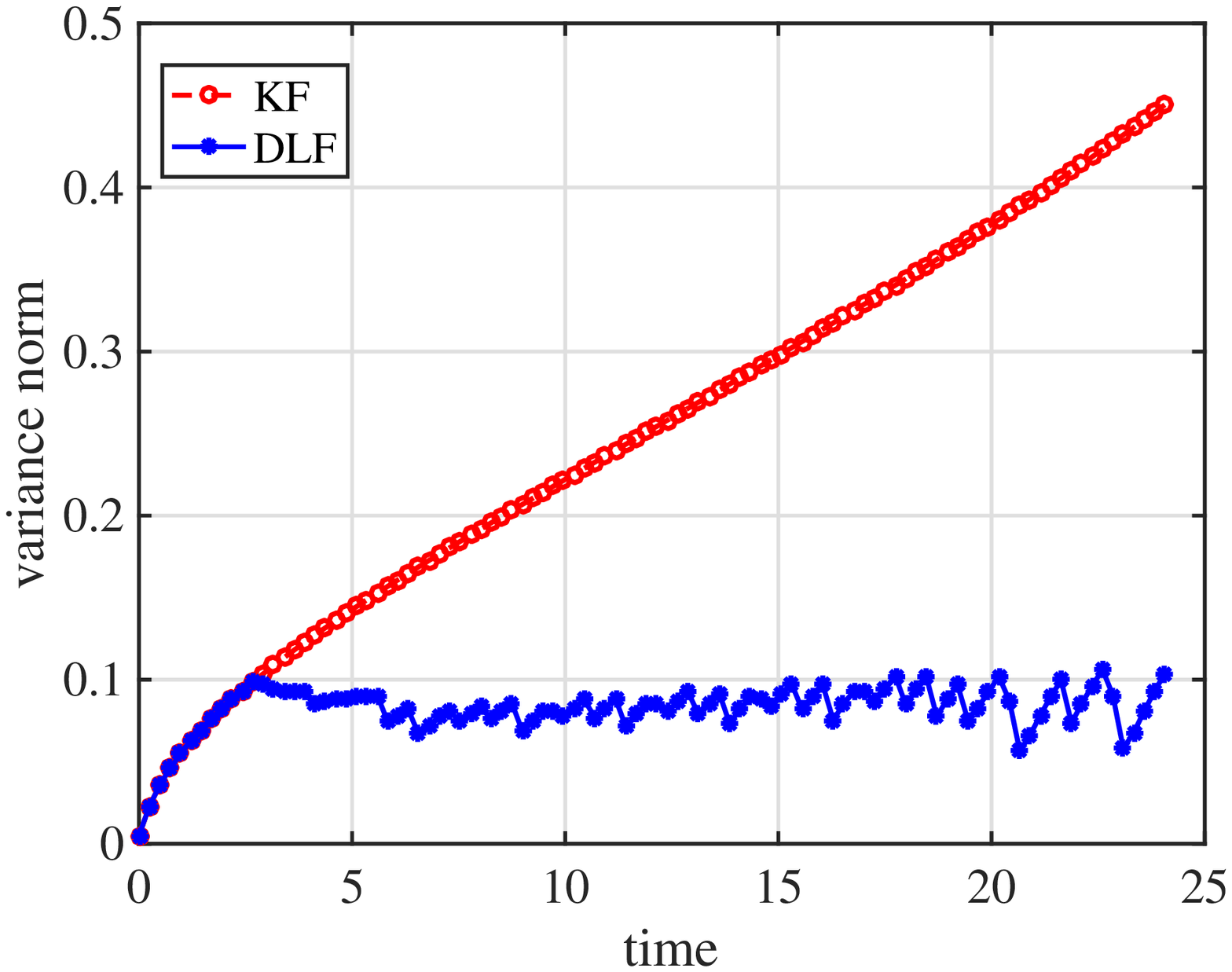}}\\
\caption{Problem II, uncertainty estimates for case (a) $\tau=1/10$, $\xi=1$, and (b) $\tau=1/10$ 
and $\xi=1/4$.}
\label{DLKFuncert}
\end{figure}

For Problem II, Figures \ref{DLKF1}-\ref{DLKFuncert} show the outcomes of  an ensemble member
of the Truth,  the Model outcome, and the KF, and DLF outcomes.   Figure \ref{DLKF1} shows the
space-time evolution corresponding 
to  $\tau=1/10$ and $\xi=1$ case; Figure \ref{DLKF2} 
to  $\tau=1$ and $\xi=1/4$, and Figure \ref{DLKF3} corresponds
to  $\tau=1/10$ and $\xi=1/4$, respectively.
The DLF and the KF are comparable
when the observation network is not sparse.  
 Figure \ref{DLKF4} shows the difference between an ensemble member of the Truth and the other three
 estimates, corresponding to the  $\tau=1/10$ and $\xi=1/4$ case shown in Figure \ref{DLKF3},  suggesting
a DLF estimate that is better than the KF estimate under sparse conditions.
The DLF and the KF have superior phase outcomes, compared to the Model. 
The time evolution of the center of mass of the estimates and an ensemble member of Truth appear in Figure \ref{DLKFphase}, for the parameters considered in Figures \ref{DLKF1}-\ref{DLKF4}.
The phase of the DLF is consistent with the Truth even when the observation system is sparse.
Although the mean estimate provided by KF is not appreciably different from the DLF for
the case when $\tau = 1$ and $\xi=1/4$, the DLF has smaller uncertainties, as evidenced in 
Figure \ref{DLKFuncert}a.
In Figure \ref{DLKFuncert}b, the uncertainties of the DLF and the KF estimates are compared. The KF
uncertainty is running away.

\section{Conclusions}
\label{discussion}

A data assimilation, specifically tailored to advectively-dominated evolution partial different equations is proposed. When compared to a more traditional filtering scheme, the method is shown to improve estimated moments of the posterior distribution, composed of a likelihood informed by noisy observations, and a prior informed by noisy model outcomes, in cases where the observation network is sparse.
 
In advection-dominated dynamics, such as is the case in models of wave-like phenomena, dramatic improvements in forecasting can be made if the phase of the wave and the features of the field are well  captured. Phase errors  arise from  epistemic errors and aleatoric errors in the phase speed. Both can affect features of the solution. 
Even if the physics are well captured by the model, the epistemic error in discrete models arises from truncation errors, which are usually deterministic. Aleatoric errors come from   uncertainties in the wave speeds. These can also affect features of the wave, compounding errors  in the wave due to amplitude uncertainties.
The focus here was on the aleatoric error, since truncation errors are practically handled by making choices on the discretization of the 
 physics  model, rather than by filtering. (Obviously, epistemic and aleatoric errors also arise in the boundary conditions and in the
initial conditions. Epistemic errors in the boundary/initial data can lead to bias errors in the estimate as well. Consideration of this source
of error is not discussed in this study, however, the aleatoric initial condition error was accounted for, in a manner that is already
standard in linear filtering methods. \cite{gilling}, for example,  considers  the implication of some forms of epistemic error and their impact on filter estimates).

Data assimilation can  lead to significant improvements in dynamic estimates, when observations with low uncertainty are blended. However, it is often the case in spatially-extended problems in geoscience that the observation system is sparse. This is
particularly challenging in problems with 3 space dimensions and time ({\it e.g.}, climate and weather models, solid mechanics, 
hydrology, etc). Moreover, with the higher resolutions afforded by larger computers the number of degrees of freedom in general 
is expected to grow by $N^3$, a rate of growth that is impractical in  observation systems. The trend is instead to decrease the 
uncertainty in the observation system. Sparse and low uncertainty observations pose a variety of challenges to filtering schemes, 
even if the biases from models and observations are properly accounted for; namely, the challenges of sampling very improbable events, 
and in hyperbolic and advection-dominated problems, in introducing phase errors in the estimates.

The proposed filtering technique, the {\it Dynamic Likelihood Filter} (DLF), is formulated  to address phase and feature errors in the estimates in advection-dominated problems, when the observation system is sparse and of low uncertainty. It
 recasts the Bayesian statement that leads to a posterior distribution of the state variable, conditioned on observations, by exploiting the hyperbolic nature of the problem to propagate  measurements forward in time. This adds a computational expense. However, since it  is linear in the number of observation stations can be  argued that in practical circumstances the added computational expense  associated with  using the deterministic time integrator to update at 
 most $N$ observations and their uncertainties is a reasonable price to pay, especially if the number of spatial dimensions reaches 3 (and thus the number of degrees of freedom, ${\cal O} (N^3)$).

For a simple linear hyperbolic system, the DLF delivered estimates that were superior to the KF, for sparse and high quality data. 
It was also shown that the data assimilation yielded estimates that were superior than estimates obtained with model-only predictions.

It is envisioned that there would be 
 two variants of the Dynamic Likelihood Filter. When the measurements are fixed in space, which is the situation that was examined in this study,  and a Lagrangian variant, in which  the measurements are taken by a moving passive or active observational platforms. 
For the fixed observation network case there are some obvious directions to pursue in future work on dynamic likelihood filtering: in addition to generalizing to higher physical spatial dimensions, the nonlinear/non-Gaussian data assimilation problems needs to be addressed. 
The present work focused on the fixed dynamic likelihood as it applies to the simplest possible hyperbolic problem. Namely, a linear one-way wave equation with uncertainties in the wave speed.  Since the problem is linear and Gaussian, a Kalman filter scheme was 
used. The smoother can be obtained by variational methods, most naturally, by the representer method ({\it cf.}, 
\cite{mead}, and \cite{bennettbook}) for the  hyperbolic dynamics case. 
Other issues that could be  pursued in future are:  (1)  how to blend data and 
models when several observations inform a single state variable, which is in fact a challenge for many other filtering schemes; (2) in DLF, how best to interpolate data in  
the multi-analysis stage; (3)  exploring  DLF forecasts,   in which the influence of data is projected forward in time with the goal 
of  using it to improve  forecasts of future states via the Bayesian methodology.

%

\ack 
This work was supported by PEER research grant \#1123-NCTRYH, and NSF/OCE \#1434198. I wish to thank Stockholm University Rossby Fellowship Program, and MPE at Imperial College, London, for their
hospitality.

\bibliographystyle{wileyqj}
\bibliography{doubwell}

\end{document}